\documentclass{amsart}
\usepackage{fullpage}
\usepackage{graphicx}
\usepackage{amsmath,amsfonts,amssymb,latexsym,epic,eepic}
\usepackage{dsfont}
\usepackage[mathcal]{euscript}
\usepackage{multicol}
\usepackage{tikz}
\usetikzlibrary{patterns}

\newcommand{\bfn}{{\boldsymbol n}}
\newcommand{\bfN}{{\boldsymbol N}_f}

\newcommand{\bfu}{{\boldsymbol u}}

\newcommand{\gradi}{{\boldsymbol \nabla}}
\newcommand{\dive}{{\rm div}}
\newcommand{\edge}{\sigma}
\newcommand{\edges}{\mathcal{E}}

\newcommand{\mesh}{{\mathcal M}}
\newcommand{\ie}{\emph{i.e.}}

\newcommand{\sumiI}{\sum_{i\in\mathcal I}}
\newcommand{\cpi}{c_{{\rm p},i}}
\begin{document}
\title[Two models for the computation of laminar flames in dust clouds]
{Two models for the computation of laminar flames in dust clouds}

\author{Dionysios Grapsas}
\address{I2M UMR 7373, Aix-Marseille Universit\'e, CNRS, \'Ecole Centrale de Marseille.}
\email{dyonisios.grapsas@univ-amu.fr}

\author{R. Herbin}
\address{I2M UMR 7373, Aix-Marseille Universit\'e, CNRS, \'Ecole Centrale de Marseille.}
\email{raphaele.herbin@univ-amu.fr}

\author{J.C. Latch\'e}
\address{Institut de Radioprotection et de S\^uret\'e Nucl\'eaire (IRSN), Saint-Paul-lez-Durance, 13115, France.}
\email{jean-claude.latche@irsn.fr}


\begin{abstract}
We address two models for the simulation of dust clouds premixed combustion: the first one consists in usual balance equations; to derive the second one, we suppose that the solution takes the form of a travelling combustion wave and track the location of the flame brush by a phase-field-like technique.
We build a finite volume fractional step scheme for both models, which respects the natural physical bounds of the unknowns.
Then we assess the consistency of both formulations.
\end{abstract}
\maketitle
%
%
\section{Introduction}\label{sec:int}

We address in this paper two alternative models dedicated to the simulation of laminar flames in dust suspensions in a gaseous atmosphere, for which a one-dimensional representation, supposing a low Mach number flow, is sufficient.
The combustible particules are supposed to be in mechanical and thermal equilibrium with the continuous phase (or, in other words, no drift nor temperature deviation between the gas and solid phases is taken into account).
We consider two descriptions of the combustion phenomenon:
\begin{list}{--}{\itemsep=0.5ex \topsep=0.5ex \leftmargin=1.cm \labelwidth=0.3cm \labelsep=0.5cm \itemindent=0.cm}
\item the first one is obtained by collecting mass balance for the chemical species, the energy balance and the momentum balance for the mixture; the reaction term $\dot \omega$ is expressed by a closure law depending of the temperature, derived on the basis of physical arguments.
This model will be refered to in the following as the {\em primitive formulation}.
\item The second one relies on the assumption that the solution consists in a travelling reaction thin interface (the so-called flame front) separating a zone where the combustion is complete (the "burnt zone") from a zone where no combustion has yet occured (the "fresh zone").
This representation offers the possibility to reduce the problem to an explicit tracking of the front location, through the solution of a transport equation for a color function $G$ ($G\in[0,1]$, $G<0.5$ in the burnt zone, $G \geq 0.5$ in the fresh atmosphere); the reaction term is governed by the value of $G$: $\dot \omega=0$ if $G \geq 0.5$ and $\dot \omega$ is proportional to $1/\tau$ otherwise, where $\tau$ is a time-scale closely correlated to the flame front thickness.
In the rest of this paper, we will call this model the {\em flame velocity formulation}.
\end{list}
The first option is standard for the computation of laminar flames.
Variants of the flame velocity formulation are often chosen to compute turbulent deflagrations in industrial applications \cite{pet-00-tur,lip-02-tur}, as in nuclear safety studies performed at the French Institut de Radioprotection et de S\^uret\'e Nucl\'eaire (IRSN).
Indeed, this latter model seems easier to solve, in the sense that stable segregated algorithms may be designed for its purpose; in addition, the flame brush incorporates structures which are very small compared to the system scales, and the flame velocity approach allows an upscaling of this complex physical phenomenon through a single parameter (the turbulent flame velocity) which may be inferred from experimental data. 
 
\medskip
A finite volume fractional step numerical scheme was developed in \cite{ami-16-mod} for the solution of the system of primitive equations; we shortly describe it here and review its main properties.
The aim of this paper is then to assess the accuracy of the switch from the first model to the second one: first, we check that the solution to the primitive formulation, in conditions representative of the target physical reality, is indeed a flame front propagating through the medium; then, we compare such a solution with the one obtained with the corresponding flame velocity model.
%
%
\section{The primitive formulation}\label{sec:mod}

\subsection{The governing equations}

The flow is supposed to be governed by the balance equations modelling a variable density flow in the asymptotic limit of low Mach number flows, namely the mass balance of the chemical species and of the mixture, the enthalpy balance, and the momentum balance equations.
For a one-dimensional flow in such a quasi-incompressible situation, the role played by the mass and momentum balance equations is quite different than in the multi-dimensional case: the velocity may be seen as the solution of the mass balance equation, and the momentum balance yields the dynamic pressure.
Since this latter unknown does not appear in the other equations, its computation is of poor interest, and the momentum balance equation may be disregarded.

\medskip
Except for this aspect, equations in this section are written in the usual multi-dimensional form.
The computational domain is denoted by $\Omega$, and its boundary $\partial \Omega$ is supposed to be split into an inflow part $\partial \Omega_I$ (where the flow enters the domain, \ie\ $\bfu \cdot \bfn_{\partial \Omega} < 0$, with $\bfu$ the flow velocity and $\bfn_{\partial \Omega}$ the normal vector to $\partial \Omega$ outward $\Omega$) and an outflow one $\partial \Omega_O$ (where the flow leaves the domain, \ie\ $\bfu \cdot \bfn_{\partial \Omega} \geq 0$) of positive $(d-1)$-measure.
The problem is posed over the time interval $(0,T)$.

\bigskip
{\bf Mass balance equations --}
The mass balance reads:
\begin{equation}\label{eq:mass}
\partial_t \rho+\dive(\rho \bfu) =0,
\end{equation}
where $\rho$ stands for the fluid density.
This equation must be complemented by an initial condition and a boundary condition on $\partial \Omega_I$ for the density; both functions are obtained by the data of the temperature and flow composition, thanks to the equation of state (see below).

\medskip
Only four chemical species are supposed to be present in the flow, namely the dust, or fuel (denoted by $F$), the oxydant ($O$), the product ($P$) of the reaction, and a neutral gas ($N$).
A one-step irreversible total chemical reaction is considered:
\[
\nu_F F + \nu_O O + N \rightarrow \nu_P P + N,
\]
where $\nu_F$, $\nu_O$ and $\nu_P$ are the molar stoichiometric coefficients of the reaction.
Chemical species other than the fuel are supposed to be gases.
The system of the mass balance equations for the chemical species reads:
\begin{equation}\label{eq:y}
\partial_t(\rho y_i)+\dive(\rho y_i \bfu) = \dot \omega_i, \qquad \mbox{for } 1 \leq i \leq N_s,
\end{equation}
where $y_i$ and $\dot \omega_i$ stand respectively for the mass fraction and the reaction rate of the species $i$.
The number of species is denoted by $N_s$, with, by assumption, $N_s=4$, and we indifferently use the notation $(y_i)_{1\leq i \leq N_s}$ or $y_F$, $y_0$, $y_P$ and $y_N$ for the fuel, oxydant, product and neutral gas mass fractions, respectively.
To simplify the exposition, the mass diffusion fluxes have been supposed in the set \eqref{eq:y} of equations to vanish.
This system must be complemented by initial and a Dirichlet boundary conditions for $(y_i)_{1\leq i \leq N_s}$ on the inflow part of the domain boundary $\partial \Omega_I$.
The prescribed values of the mass fractions at the initial time and on the inflow boundary lie in the interval $[0,1]$.
The reaction rate of each chemical species may be written as:
\[
\dot \omega_F = -\nu_F W_F\, \dot \omega,\quad \dot \omega_O = -\nu_O W_0\, \dot \omega,
\quad \dot \omega_P = \nu_P W_P\, \dot \omega
\quad \mbox{and } \dot \omega_N = 0,
\]
where $W_F$, $W_0$ and $W_P$ stand for the molar masses of the fuel, oxydant and product respectively, and $\dot \omega$ is a non-negative reaction rate, which is supposed to vanish when either $y_F=0$ or $y_O=0$.
Since $\nu_F W_F + \nu_O W_0 = \nu_P W_P$, we have $\dot \omega_F + \dot \omega_O + \dot \omega_P =0$.

\bigskip
{\bf Energy balance --}
In the low Mach number approximation, the total enthalpy balance reads:
\begin{equation}\label{eq:T}
\sumiI \cpi\ \Bigl[ \partial_t (\rho y_i \theta) + \dive(\rho y_i \theta \bfu) \Bigr]
-\dive(\lambda\gradi \theta)=\dot \omega_\theta,
\quad \dot \omega_\theta = -\sumiI \Delta h_{f,i}^0 \dot{\omega}_i.
\end{equation}
where $\theta$ stands for the temperature, $\cpi$ for the specific heat of the species $i$ (supposed to be constant), $\Delta h_{f,i}^0$ for the formation enthalpy at $0^\circ$K and $\lambda$ the thermal conductivity.
This equation is complemented by a total flux boundary condition at the inlet boundary, and we suppose that the diffusion flux vanishes at the outlet boundary.

\bigskip
{\bf Equation of state --}
We suppose that the gas phase is a mixture of perfect gases and that the density $\rho_F$ of the solid phase (\ie\ of the fuel) is constant, so:
\begin{equation}\label{eq:eos}
\rho = \varrho \bigl(\theta,(y_i)_{1\leq 1 \leq N_s}\bigr)
= \frac 1 {\displaystyle \frac{R \theta}{P_{th}}\sum_{i=O,P,N} \frac{y_i}{W_i} + \frac{y_F}{\rho_F}},
\end{equation}
where $R=8.31451\,\mbox{JK}^{-1}\mbox{mol}^{-1}$ stands for the perfect gases constant.
Since the computational domain is supposed not to be closed, the so-called thermodynamic pressure $P_{th}$ is constant in time and space, and given by the initial state.
%
%
\subsection{The numerical scheme}\label{sec:scheme}
For the solution of the equations of the model, we define the variable $z$ as follows:
\[
z = \frac{s\,y_F+1-y_O}{1+s}, \quad \mbox{with } s=\frac{\nu_O W_O}{\nu_F W_F}.
\]
Note that, combining the mass balance equation for the fuel and the oxydant, the variable $z$ satisfies an homogeneous equation; for this reason, we replace the oxydant mass balance equation by the balance equation for $z$ (since, given the values of $z$ and $y_F$, we may deduce $y_O$).

\medskip
Let us consider a partition $0=t_0 < t_1 <\ldots < t_N=T$ of the time interval $(0,T)$, which we suppose uniform.
Let $\delta t=t_{n+1}-t_n$ for $n=0,1,\ldots,N-1$ be the constant time step.
We suppose that the interval $\Omega$ is split into a family of control volumes (sub-intervals of $\Omega$) which realizes a partition of $\Omega$; we denote these control volumes by $(K)_{K\in\mesh}$.
The scalar unknowns, \ie\ the density, mass fractions and temperature, are associated to the control volumes, and the corresponding unknowns read $\rho_K^n$, $(y_i)_K^n$, $z_K^n$ and $\theta_K^n$ for $K \in \mesh$, $0 \leq n \leq N$ and $i\in\mathcal I$.
The velocity is discretized at the faces of the mesh, which we denote by $(\edge)_{\edge\in \edges}$, so the corresponding unknowns are $u^n_\edge$ for $\edge \in \edges$ and $0 \leq n \leq N$.
We implement a fractional-step algorithm, which consists in four steps, in order to calculate recursively the unknowns $(y_i)^{n+1}_{i\in\mathcal I}$, $z^{n+1}$, $\theta^{n+1}$, $\rho^{n+1}$ and $u^{n+1}$ for $0\leq n < N$:
\begin{subequations}\label{eq:scheme_P}
\begin{align}
\nonumber &
\mbox{{\bf Chemistry step} -- Solve for $(y_N, z, y_F, y_P)^{n+1}$:}
\\ \nonumber & \quad
\forall K \in \mesh,
\\[0.5ex] \label{eq:sch-yN} & \hspace{3ex}
\dfrac 1 {\delta t} \bigl[\rho^n_K (y_N)_K^{n+1}-\rho^{n-1}_K (y_N)_K^n \bigr]
+ \dive \Bigl[\rho^n y_N^{n+1} \bfu^n \bigr]_K
=0,
\\[0.5ex] \label{eq:sch-z} & \hspace{3ex}
\dfrac 1 {\delta t} \bigl[\rho^n_K z_K^{n+1}-\rho^{n-1}_K z_K^n \bigr]
+ \dive \bigl[\rho^n z^{n+1} \bfu^n \bigr]_K
= 0,
\\[0.5ex] \label{eq:sch-yF} & \hspace{3ex}
\dfrac 1 {\delta t} \bigl[\rho^n_K (y_F)_K^{n+1}-\rho^{n-1}_K (y_F)_K^n \bigr]
+ \dive \bigl[\rho^n y_F^{n+1} \bfu^n \bigr]_K
=(\dot\omega_F)_K^{n+1},
\\[0.5ex] \label{eq:sch-yP} & \hspace{3ex}
(y_F)_K^{n+1} + (y_O)_K^{n+1} + (y_N)_K^{n+1} + (y_P)_K^{n+1} = 1.
\displaybreak[1] \\[3ex] \nonumber &
\mbox{{\bf Energy balance} -- Solve for $\theta^{n+1}$:}
\\ \nonumber & \quad
\forall K \in \mesh,
\\[0.5ex] \label{eq:sch-e} & \hspace{3ex}
\begin{matrix}\displaystyle
\sumiI \cpi \Bigl[
\dfrac 1 {\delta t} \bigl[\rho^n_K (y_i)_K^{n+1} \theta_K^{n+1}-\rho^{n-1}_K (y_i)_K^n \theta_K^n \bigr]\Bigr]
+ \dive \bigl[\rho^n y_i^{n+1} \theta^{n+1} \bfu^n \bigr]_K
\hspace{10ex} \\ \hfill
-\dive(\lambda\gradi \theta^{n+1})_K
=(\dot\omega_\theta)^{n+1}_K.
\end{matrix}
\displaybreak[1] \\[3ex] \label{eq:sch-rho} &
\mbox{{\bf Equation of state} -- $\rho_K^{n+1} =\varrho \bigl(\theta_K^{n+1},((y_i)_K^{n+1})_{1\leq 1 \leq N_s}\bigr)$, for $K \in \mesh$.}
\displaybreak[1] \\[3ex] \nonumber &
\mbox{{\bf Mass balance} -- Solve for $\bfu^{n+1}$:}
\\ \label{eq:sch-mass} & \quad
\forall K \in \mesh,\qquad
\dfrac 1 {\delta t} \bigl[\rho_K^{n+1}-\rho_K^n \bigr]
+ \dive \bigl[\rho^{n+1} \bfu^{n+1} \bigr]_K
=0.
\end{align}
\end{subequations}

\medskip
The discrete operators appearing in these relations are approximated by finite-volume techniques.
Thanks to a careful definition of the convection fluxes, derived to fulfill the conditions introduced in \cite{lar-91-how} to obtain a maximum-principle-preserving convection operators, the scheme is proven in \cite{ami-16-mod} to preserve the physical bounds of the unknowns: mass fractions in the interval $[0,1]$, positivity of the temperature and so of the density.
%
%
\section{A model based on an explicit tracking of the flame front}\label{sec:mod:level_set}

\subsection{The governing equations}

The physical model addressed in this section is based on an explicit computation of the flame brush location, by a phase-field-like technique.
The "color function" is called $G$ in this context; its transport equation is referred to as the $G$-equation, and reads:
\begin{equation}\label{eq:G}
\partial_t(\rho G) + \dive (\rho G \bfu) + \rho_u u_f |\gradi G| =0,
\end{equation}
Initial conditions are $G=0$ at the location where the flame starts and $G=1$ elsewhere.
The quantity $\rho_u$ is a constant density, which, from a physical point of view, stands for a characteristic value for the unburnt gases density, and $u_f$ is the flame brush velocity.
The reactive term $\dot \omega$ is given by:
\begin{equation}\label{eq:wG}
\dot \omega= \frac{u_f}{\delta}\ \eta(y_F,y_O)\ (G-0.5)^-,\quad \eta(y_F,y_O)=\min (\frac{y_F}{\nu_F W_F}, \frac{y_O}{\nu_O W_O}),
\end{equation}
where $\delta$ is a quantity homogeneous to a length scale, which governs the thickness of the reaction zone.

\medskip
The flame velocity model consists of Equation \eqref{eq:G}, of the mixture mass balance equation \eqref{eq:mass}, the chemical species mass balance equations \eqref{eq:y} (with the modified expression \eqref{eq:wG} for the chemical reaction term $\dot \omega$) and of the energy balance \eqref{eq:T}.
Note that, under some assumptions which are usually not valid in industrial applications, this model may be simplified: for instance, in perfectly premixed situations (\ie\ constant in space initial data for the chemical mass fractions and the temperature) and supposing an infinitely fast chemical reaction (\ie, in the present formalism, making $\delta$ tend to zero), the variable $G$ may be indentified to a progress variable and all the other unknowns (\ie\ the mass fractions and the temperature) may be deduced from $G$ through an algebraic relation.

\subsection{Numerical scheme}

The $G$ function is discretized on the primal mesh, so the discrete unknowns are $G^n_K$, for $K\in\mesh$ and $0 \leq n \leq N$.
The numerical algorithm differs from the scheme for the primitive formulation, \ie\ System \eqref{eq:scheme_P}, by the insertion, as a first step, of a discrete counterpart to Equation \eqref{eq:G}:
\begin{subequations}
\begin{align}
\nonumber &
\mbox{{\bf flame brush transport step} -- Solve for $G^{n+1}$:}
\\ \nonumber & \quad
\forall K \in \mesh,\qquad
\dfrac 1 {\delta t} \bigl[\rho^n_K G_K^{n+1}-\rho^{n-1}_K G_K^n \bigr]
+ \dive \Bigl[\rho^n G^{n+1} \bfu^n \bigr]_K
+ \rho_u u_f |\gradi G|^{n+1}_K
=0.
\end{align}\end{subequations}
For the discretization of the last term in this relation, we write:
\[
|\gradi G| = \frac{\gradi G}{|\gradi G|} \cdot \gradi G, \qquad \mbox{so }\quad
|\gradi G|^{n+1}_K = (\bfN^n \cdot \gradi G^{n+1})_K,
\]
where $\bfN$ is an approximation of the advection field $\gradi G/|\gradi G|$ and we use an upwind finite volume formulation of the transport operator, \ie\ the formulation obtained by writing $\bfN \cdot \gradi G = \dive(G\,\bfN) - G\,\dive \bfN$ and using an upwind finite volume (first or second order) discretization of the convection operator.
For the present solver, the convection operators are discretized by an explicit MUSCL-like technique \cite{pia-13-for}.
%
%
\section{Results}\label{sec:results}

\begin{figure}[tb]
\begin{center}
\includegraphics[width=0.8\textwidth]{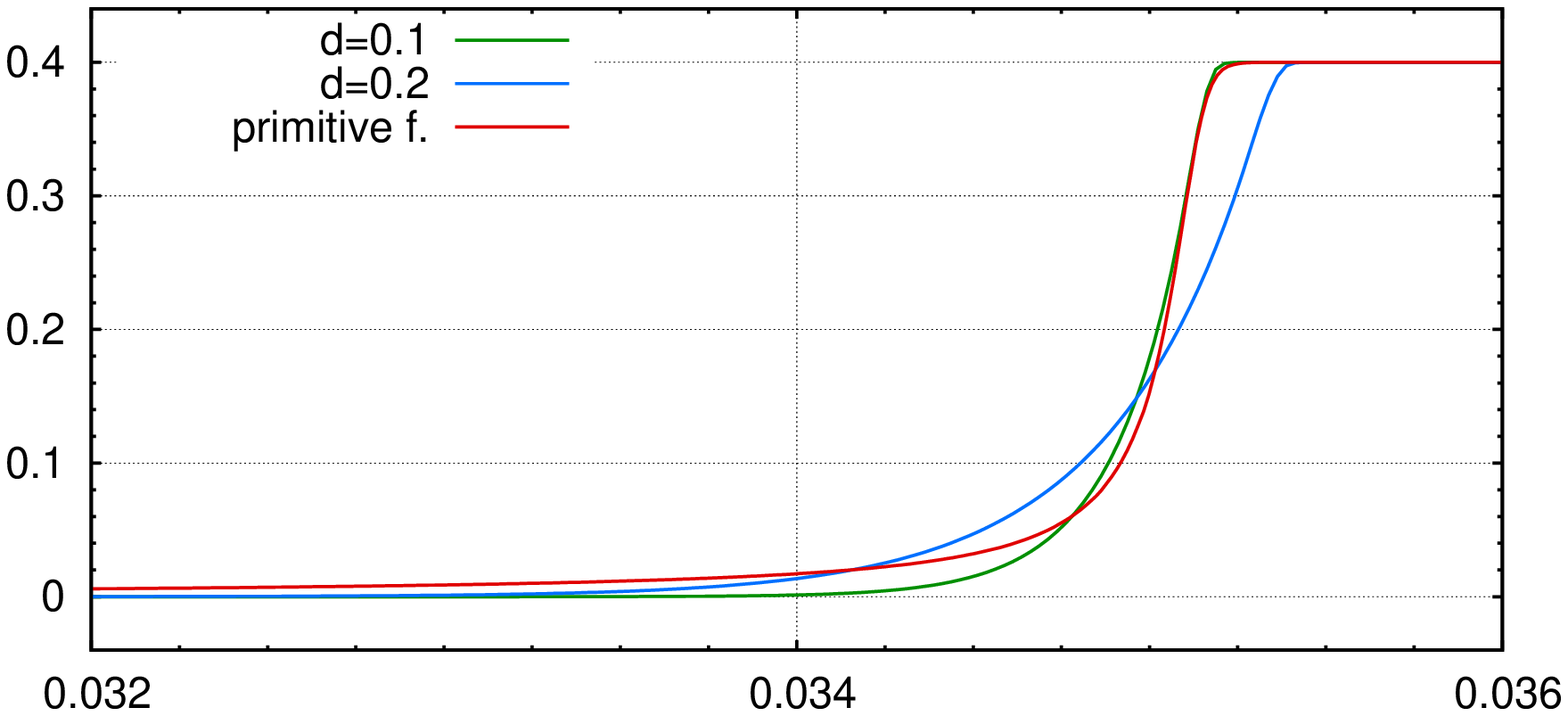}\\
\includegraphics[width=0.8\textwidth]{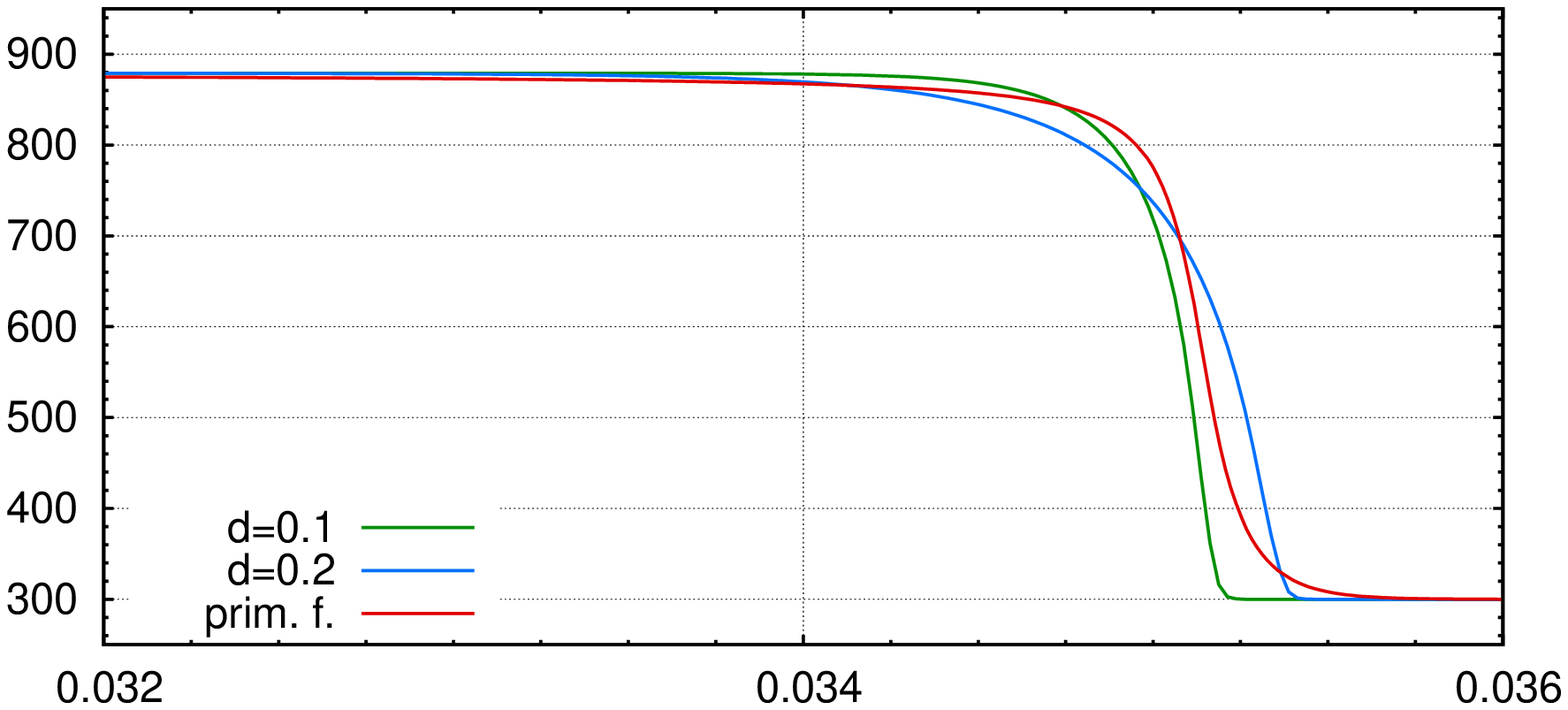}
\end{center}
\caption{Fuel mass fraction (top) and temperature (bottom) travelling profiles obtained with the primitive formulation of the equations (red) and with the flame velocity model, with $\delta=0.1$\,mm (green) and $\delta=0.2$\,mm (blue).}
\label{fig:profiles} \end{figure}

Computations presented in this section are performed with MATLAB for the primitive formulation and by the open-source CALIF$^3$S software developped at IRSN \cite{califs} for the flame-velocity model.

\smallskip
Data is chosen in order to allow to check the scheme properties (\ie\ to avoid unrealistic simplifications, as, for instance, a same specific heat diffusion coefficient for all the chemical species), and to be in the range of practical applications.
The mixture is initially at rest, homogeneous and with a uniform temperature:
\[
 (y_F)_0=(y_O)_0=0.4,\quad (y_N)_0=0.2,\quad (y_P)_0=0,\quad \theta_0=300^\circ K.
\]
In the primitive formulation, the reaction rate follows an Arrhenius law:
\[
\dot\omega_K = 10^4\,y_F\, y_O\ e^{-900/\theta}.
\]
The molar masses of the chemical species are considered to be equal to $20$\,g mol$^{-1}$ for all the species, so the combustion reaction reads $F+O+N \longrightarrow 2P+N$, and the initial atmosphere composition is stoichiometric.
The temperature diffusion coefficient is $\lambda=0.005$, the specific heat coefficients (J Kg$^{-1}$\ K$^{-1}$) are $c_{{\rm p},N}=3.\,10^3$, $c_{{\rm p},F}=1.\,10^3$, $c_{{\rm p},O}=2.\,10^3$ and $c_{{\rm p},P}=4.\,10^3$ and the formation enthalpies (J Kg$^{-1}$) are $\Delta h_{f,N}^0=3.\,10^6$, $\Delta h_{f,F}^0=1.\,10^6$, $\Delta h_{f,O}^0=-2.\,10^6$ and $\Delta h_{f,P}^0=-4.\,10^6$ (so the reaction is exothermic).
The fuel density is equal to $100$\, Kg m$^{-3}$.
Ignition is obtained in the primitive formulation by making $\dot \omega$ depend in a very thin zone on a fictitious elevated temperature, to trigger a reaction at the initial time.
In the flame velocity model, $G$ is imposed to zero in the same zone (while $G=1$ elsewhere).
Since the inflamation zone is very thin, the consequent initial burst is not too violent.

\smallskip
First of all, we observe that, for the primitive formulation, the solution tends after a transition period to a travelling combustion wave separating a fresh (or unburnt) zone from a burnt zone, where $y_F=y_O=0$, $y_P=0.8$ and the temperature is equal to the adiabatic combustion temperature.
By construction of the scheme, the neutral gas mass fraction $y_N$ and the reduced variable $z$ are kept constant in time and space and equal to their initial value.
Since the profile in the interface does not vary in space and time up to a translation velocity $u_p$ (the velocity of the flame brush), we may write the jump conditions for the mixture mass balance equation, to obtain:
\[
(\rho_u -\rho_b)\ u_p = \rho_u u_u -\rho_b u_b,
\]
where $\rho_b$ and $u_b$ (resp. $\rho_u$ and $u_u$) stand for the constant density and velocity in the burnt (resp. unburnt) zone.
Thanks to symmetry conditions (due to the fact that the combustion takes place in an atmosphere initially at rest), $u_b=0$ and we deduce from the previous relation that the flame velocity is given by:
$
u_f=u_p -u_u = u_u \ \rho_b/(\rho_u-\rho_b).
$
The obtained value is injected in the flame velocity model, and we choose the length $\delta_f$ to fit as closely as possible the travelling profiles of the unknowns.
Results for the fuel mass fraction and the temperature are given on Figure \ref{fig:profiles}.
We observe that, as expected, the thickness of the combustion zone is scaled by $\delta_f$ and that a reasonable agreement is obtained with $\delta_f=0.1$\,mm.
%
%
\bibliographystyle{plain}
\bibliography{dust}
\end{document}